\documentclass{article}
\usepackage{amssymb}
\usepackage{epsfig}

\usepackage{color}

\newcommand{\eq}{\begin{equation}}
\newcommand{\en}{\end{equation}}
\newcommand{\re}[1]{\mbox{(\ref{#1})}}

\newtheorem{Theorem}{Theorem}
\newtheorem{theorem}[Theorem]{Theorem}
\newtheorem{lemma}[Theorem]{Lemma}
\newtheorem{corollary}[Theorem]{Corollary}
\newtheorem{construction}[Theorem]{Construction}
\newtheorem{proposition}[Theorem]{Proposition}
\newtheorem{example}[Theorem]{Example}
\newtheorem{defn}[Theorem]{Definition}

\newtheorem{question}[Theorem]{Question}
\newtheorem{conjecture}[Theorem]{Conjecture}
\newtheorem{condition}[Theorem]{Condition}
\newtheorem{remark}[Theorem]{Remark}
\newtheorem{problem}[Theorem]{Problem}
\def\proof{\noindent{\bf Proof.\ \ }}

\def\endpf{\hfill $\square$ \vskip .25in}

\def\cmaj {concave majorant}

\newcommand {\ovC} {\overline{C}}

\newcommand {\hx} {\hat{x}}

\newcommand {\hu} {\hat{u}}
\newcommand {\pu} {u'}
\newcommand {\hhat} {\hat{t}}
\newcommand {\pt} {t'}
\newcommand {\hy} {\hat{y}}

\newcommand{\IP}{\mathbb{P}}
\newcommand{\IE}{\mathbb{E}}

\newcommand{\ints}{\mathbb{Z}}

\newcommand{\te}{\rightarrow}
\newcommand{\ed}{\mbox{$ \ \stackrel{d}{=}$ }}

\newcommand{\eps}{\varepsilon}

\newcommand{\lb}[1]{\label{#1}}

\newenvironment{thm}[1]{\begin{theorem}\label{#1}}{\end{theorem}}

\newenvironment{lmm}[1]{\begin{lemma}\label{#1}}{\end{lemma}}

\newenvironment{crl}[1]{\begin{corollary}\protect\label{#1}}{\end{corollary}}

\newenvironment{prp}[1]{\begin{proposition}\protect\label{#1}}{\end{proposition}}

\newcommand{\br}{B^{\mbox{$\scriptstyle{\rm br}$}}}

\def\proof{\noindent{\bf Proof.\ \ }}

\def\endpf{$\Box$}

\begin{document}

\title{The distribution of the maximal difference between Brownian bridge and its concave majorant}

\author{Fadoua Balabdaoui\thanks{Universit\'e  Paris-Dauphine, Place du Mar\'echal de Lattre de Tassigny, 75775 Paris CEDEX 16, France} \ and Jim Pitman\thanks{University of California, Berkeley; research supported in part by N.S.F.\ Grant DMS-0806118}}

\maketitle

\begin{abstract} 

We provide a representation of the maximal difference between a standard Brownian bridge and its \cmaj \ on the unit interval, from which we deduce expressions for the distribution and density functions and moments of this difference.
This maximal difference has an application in nonparametric statistics where it arises in testing monotonicity of a density or regression curve.


%
\end{abstract}

\section{Introduction}

\begin{figure}[!h]
\centering
\centerline{\epsfig{file = 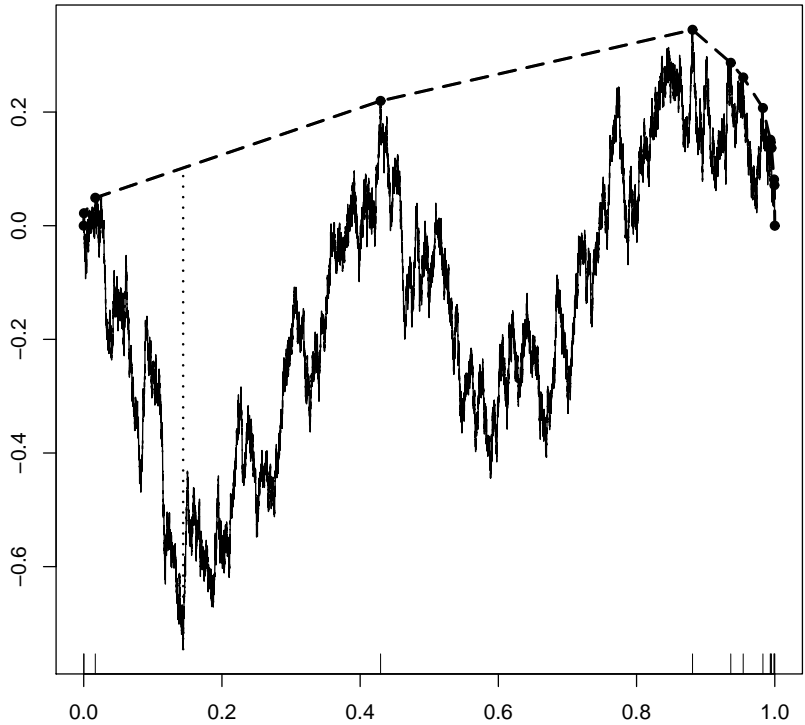, width = 12cm}}
\vspace*{-0.5cm}
\caption{Plot of a Brownian bridge and its concave majorant. The bullets depict the vertices, whose times are indicated by the tick marks on the x-axis. The length of the dotted vertical segment gives the maximum difference between the Brownian bridge and the \cmaj. The Haar approximation was used to generate the Brownian bridge on a discrete partition of $[0,1]$ with a mesh equal to $2^{-15}$.}
\label{fig: LCMBB}
\end{figure}

Motivated by applications to the theory of nonparametric statistics, indicated later in this introduction, 
we provide a useful representation of the maximal difference
$$
M := \sup_{u \in [0, 1]} \left([\ovC_{[0,1]} B](u) - B(u)\right) 
$$
where $(B(u), 0 \le u \le 1)$ is a Brownian motion, and $\ovC_I f$ denotes the (least) \cmaj \ of a function $f$ defined on an interval $I$. 
See Figure \ref{fig: LCMBB}.

Our representation of $M$ is presented in the following theorem,
in terms of the distribution of $M_3$, the maximum of a standard Brownian excursion,
which can be represented also as the maximum of a three dimensional Bessel bridge, or as
\eq
\label{m3diff}
M_3 \ed \sup_{u \in [0,1]} \br(u) - \inf_{u \in [0,1]} \br(u) 
\en
where  $\br$ denotes the {\em standard Brownian bridge} obained by conditioning $B$ on $B(1) = 0$. 
See \cite{MR0402955}, \cite{bpy99z} for background and further information about the distribution of $M_3$, which has been extensively studied.
\begin{thm}{Rep1}
The distribution of the maximal difference $M$ between a Brownian motion $B$
and its concave majorant is determined by the identity in distribution
\eq
\label{mainid}
M \ed \max_{j}  \sqrt{L_j} M_{3,j}
\en
where 
\begin{itemize}
\item
$(L_j, j = 1,2 \ldots)$ is the {\em uniform stick-breaking process}
$$
L_1:=  W_1, ~~L_2:= (1-W_1) W_2, ~~L_3:= (1-W_1) (1 - W_2) W_3, \ldots 
$$
derived from a sequence $W_1, W_2, \ldots $ of independent uniform $(0,1)$ variables;
\item 
$(M_{3,j}, j = 1,2 \ldots)$ is a sequence of independent random variables, each with the distribution of $M_3$;
\item  the two sequences $(W_j)$ and $(M_{3,j})$ are independent.
\end{itemize}
Moreover, $M$ is independent of $B(1)$, so the identity in distribution {\rm \re{mainid}} holds also 
if the Brownian motion $B$ on $[0,1]$ is replaced by the Brownian bridge obtained by
conditioning $B$ on $B(1) = b$, for arbitrary real $b$.
\end{thm}

To be more explicit, following 
Groeneboom \cite{MR714964},
we observe that the vertices of the concave majorant $([\ovC_{[0,1]} B](u), 0 \le u \le 1)$ partition
$[0,1]$ into a countable collection of subintervals, with accumulation of vertices at both $0$ and $1$, but with only a finite number of vertices almost surely in
$(\epsilon, 1 -\epsilon)$ for each $\eps >0$. The sequence $(L_j)$ represents the lengths of these maximal subintervals,
over each of which the \cmaj \ has a segment with particular slope. These $L_j$ are arranged in a suitable
random order, while $\sqrt{L_j} M_{3,j}$ represents the maximum value of $[\ovC_{[0,1]} B](u) - B(u)$ for $u$ in the interval of length $L_j$. Here and throughout the paper, the {\em length} of a segment of a concave majorant
refers to the length of the time interval associated with the segment, rather than the length of the segment in a
two-dimensional picture such as Figure \ref{fig: LCMBB}.
\medskip

For the proof of Theorem \ref{Rep1} we combine two different ingredients:
\begin{itemize}
\item Groeneboom's description 
\cite{MR714964}
of the concave majorant of Brownian motion on the infinite interval $[0,\infty)$;
\item Suidan's description \cite{suidan01} 
of the joint law of ranked lengths of intervals in the  partition of the time interval $[0,1]$ generated by vertices of $\ovC_{[0,1]} B$.
\end{itemize}
In principle, the second of these ingredients must be derivable from the first. Following Groeneboom, we use Doob's transformation to map Brownian motion 
on $[0,\infty)$ to Brownian bridge on $[0,1]$, and this mapping determines the law of the interval partition of $[0,1]$ derived from $\ovC_{[0,1]} \br$. While we provide some details of this in Section \ref{complements}, we are 
unable to fully derive the stick-breaking representation of interval lengths this way.  Still, by developing the results of Groeneboom and Suidan, and by exploiting the fact that a uniform stick-breaking process is invariant under a size-biased random permutation
 (see \cite{MR867196},\ 
\cite{MR1387889}, and 
\cite{MR1434129}) 
we see that the $L_j$ appearing in Theorem \ref{Rep1} can be constructed as follows:
\begin{crl}{uniform}
Let $U_1, U_2, \ldots$ be a sequence of independent uniform $(0,1)$ variables, independent of the Brownian motion $B$, and let $L_{1}, L_{2}, \ldots$ be the size-biased random 
permutation of lengths of segments of $\ovC_{[0,1]} B$ defined by:
\begin{itemize}
\item $L_1$ is the length of the interval containing $U_1$, 
\item $L_2$ is the length of the interval containing the first point $U_i$ that is not in this interval,
\item $L_3$ is the length of the interval containing the next point $U_j$ that is not in either of the first two intervals, 
\end{itemize}
and so on.  Then 
$$
L_1, L_2/(1 - L_1), L_3/(1 - L_1 - L_2), \ldots
$$
is a sequence of independent uniform $(0,1)$ random variables, and this sequence is independent of $B(1)$. 
\end{crl}
In particular, 
\begin{itemize}
\item
the length $L_1$ of the segment of $\ovC_{[0,1]}\br$ covering a uniform random number $U_1$ is itself uniformly 
distributed in $[0,1]$.
\end{itemize}
We show in an appendix how this fact can be verified from 
Groeneboom's  \cite{MR714964} 
joint density for the location of vertices of the concave majorant of
Brownian motion on $[0,\infty)$. But the computations are difficult and we are currently unable to extend this method to obtain the assertions about $L_j$ for $j \ge 2$ from Groeneboom's results. Rather, we 
rely completely on Suidan's approach for this part of the argument, which is essential for our proof of Theorem \ref{Rep1}.


%
%
%
%

To motivate our study of the distribution of $M$, we recall that in testing whether a density or a regression function on $[0,1]$ is decreasing, the supremum distance between the empirical estimator of the function and its concave majorant is used as a test statistic; see Kulikov and Lopuha\"a \cite{MR2211085} and Durot  \cite{MR1996191}. 
From \cite{MR2391249} it follows that the the supremum distance between the empirical distribution function and its \cmaj \ attains its
maximum at the uniform density on $[0,1]$.
Moreover, it is known that this statistic, when multiplied by $\sqrt{n}$, converges in distribution
to the maximal difference $M$ between a Brownian motion $B$ on $[0,1]$ and its \cmaj. For the regression problem, the distribution 
of $M$ appears again as the limit of a similar scaling of the supremum distance between the cumulative sum diagram and its \cmaj \ if the true decreasing regression is 
constant, which is also known to be the least favorable regression function for the testing problem in question \cite{MR1996191}. 
Durot \cite{MR1996191} established continuity of the distribution of $M$. 
It is an immediate corollary of Theorem \ref{Rep1} that $M$ in fact has a density, and we provide formulas for this density and for moments of $M$
in Section \ref{distM}.
%
Additionally we give an alternative characterization of $M$ based on the inverse of the Laplace transform of a function involving modified Bessel functions of the second kind.  
This other characterization suggests a way of calculating the quantiles of $M$ at any desired precision using some appropriate approximation method for the inverse 
of a Laplace transform, but we do not pursue this here.



  

\section{The concave majorant of Brownian motion on a finite interval}


We need to show that conditionally given $\ovC_{[0,1]} B$ the difference process 
$\ovC_{[0,1]} B - B$  behaves like a succession of independent Brownian excursions between the vertex times of $\ovC_{[0,1]} B$.  
Groeneboom \cite{MR714964} proved
a corresponding result for a Brownian motion on $[0, \infty)$, stated as Theorem \ref{excns} below. From this result on $[0, \infty)$, we will prove a similar theorem on $[0,1]$ by using a space-time transformation. 

For a fixed $a > 0$ let
$$
Z_a : = \max_{t \ge 0} \{B(t) - a t \} 
$$ 
and $D_a$ the time at which the maximum is attained; see also \cite{brownistan}. The point $(D_a, Z_a + a D_a)$ is one vertex of $\ovC_{[0,\infty)} B$, and it follows from \cite[Theorem 2.1]{MR714964} that $(D_{1/a}, a > 0)$ is a pure jump process. Let $S_0 > S_1 > S_2 > \cdots $  denote the successive slopes of  $\ovC_{[0,\infty)} B$ \ to the right of $D_a$, and 
$S_{-1} <  S_{-2} < \cdots $  the successive slopes of $\ovC_{[0,\infty)} B$ \ to the left of $D_a$, so $S_0 < a < S_{-1}$ almost surely.   For $i \in \ints$ let $T_i$ denote the length of the interval on which the slope of the \cmaj \ is $S_i$, so that $\{V_i : = \sum_{j \le i} T_j, i \in \ints \} $ is the sequence of times of vertices of $\ovC_{[0,\infty)} B$.

\begin{thm}{excns}(Groeneboom \cite{MR714964})
Conditionally given $\ovC_{[0, \infty)} B$,
the difference process 
$(\ovC_{[0, \infty)} B(t) - B(t), t \ge 0)$,
is a succession of independent Brownian excursions of prescribed lengths $T_i = V_i - V_{i-1}, i \in \ints$ between zeros at the times $V_i, i \in \ints$.
\end{thm}


See \cite{MR714964}, \cite{MR733673}, \cite{MR770946} for various proofs of this result. 

%
%

\medskip

It is important now to distinguish clearly the restriction to $[0,1]$ of 
$\ovC_{[0,\infty)}B$, the \cmaj \ of a Brownian motion $B$ on $[0, \infty)$, 
and $\ovC_{[0,1]}B$, the \cmaj \ of a Brownian motion $B$ on $[0, 1]$. The former concave majorant has vertices accumulating only at $0$, whereas the latter,
which is the subject of the following Proposition \ref{indep}, has vertices accumulating at both $0$ and $1$. These two concave majorant processes agree on some 
random interval $[0,R]$ with $0 < R < 1$ almost surely, and then differ on $(R,1]$, where $R$ is the time of the last vertex of $\ovC_{[0,\infty)}B$ before time 1, whose distribution is determined by formula (\ref{straddt}). 
See also the remark of Groeneboom \cite[p. 1022]{MR714964}.

\begin{prp}{indep} 
Let $\{V_{i}, i \in \ints \}$  denote an indexing of the times of vertices of  $\ovC_{[0,1]} B$ so that $V_i$ is an increasing function of $i \in \ints$ with $\lim_{i \te - \infty } V_i = 0$ and $\lim_{i \te  \infty } V_i = 1$. Moroever, we assume that this indexing depends only on $\ovC_{[0,1]} B$.

\begin{itemize}
\item [(i)]
The random set of times of vertices $\{ V_{i}, i \in \ints \}$ is independent of $B(1)$.
\item [(ii)]
The difference process $(\overline{C}_{[0,1]}B(u) - B(u), 0 \le u \le 1)$ is independent of $B(1)$. 
\item [(iii)]
Consequently, the distribution of both the random set of vertex times and of the difference process is the same for
an unconditioned Brownian motion $B$ as it is for a Brownian bridge from $(0,0)$ to $(1,b)$ for every real value of $b$.
\end{itemize}
\end{prp}

\medskip

\proof 
The first two assertions follow easily from the facts that
\begin{itemize}
\item
the set of times of vertices of the \cmaj \ of a function $(f(u), 0 \le u \le 1)$ is the same for $f(u)$
as it is for $f(u) + c u$ for arbitrary real $c$;
\item
for a Brownian motion $B$, the process $(B(u) - u B(1), 0 \le u \le 1)$ is a standard Brownian bridge independent of $B(1)$.
\end{itemize}
The third assertion follows easily from the first two. 
\endpf

\begin{thm}{conc}
Let $X$ be either an unconditioned Brownian motion, or a Brownian bridge from $(0,0)$ to $(1,b)$ for some real $b$, and  $\{V_i, i \in \mathbb Z\}$ an indexing of the times of vertices of $\ovC_{[0,1]} X$ as in Proposition \ref{indep}. Then conditionally given  $V_{i} = v_i, i \in \ints$ the difference process 
$(\overline{C}_{[0,1]}X(u) - X(u), 0 \le u \le 1)$ 
is a concatenation of independent Brownian excursions of prescribed lengths $v_i - v_{i-1}$ between zeros at the times $v_i$.
\end{thm}

\medskip

\par \noindent \textbf{Proof.} \ Based on Proposition \ref{indep} (iii), it suffices to consider the case with $B$ replaced by
a standard Brownian bridge $\br$.
Then, according to Doob's transformation,
$$
\br(u) = (1-u) \widehat B\left(\frac{u}{1-u} \right), \ u \in [0,1] 
$$
where $\widehat B$ is the standard Brownian motion
$$
\widehat B(t) = (1+t) \br\left(\frac{t}{1+t}\right), \ t \ge 0.
$$
It follows that the times $V_i$ of vertices of the \cmaj \ of $\br$ on $[0,1]$ may be constructed as $V_i = T_i/(T_i + 1)$
where the $T_i$'s are the times of vertices of the \cmaj \ of $\widehat B$ on $[0,\infty)$. Moreover, the difference between  $\br$ and its \cmaj \ on $(V_{i}, V_{i+1})$ is a transformation of the difference between $\widehat B$ and its \cmaj \ on $(T_{i}, T_{i+1})$.

We will show now that this transformation maps Brownian excursions to Brownian excursions. Indeed, observe  that the transformation between $(u, \br(u)) $ and $(t, \widehat B(t))$ is the restriction to a Brownian path of the
space-time transformation
\eq
\lb{spacetime}
\mathcal T(u,x) = (t,y) :=  \left(\frac{u}{1-u}, \frac{x}{1-u} \right)
\en
where $0 \le u < 1$, $0 \le t <  \infty$, and $x$ and $y$ both range over all real numbers.

Consider now the conditioning of $\br$ on $\br(u) = x$ and $\br(\hu) = \hx$ for some $0 < u < \hu < 1$ and real numbers $x$ and $\hx$. The process
$$
X^*(v):= \frac{\br( u + v ( \hu - u )) - x - v ( \hx - x )}{\sqrt{\hu - u }}, ~~~~ 0 \le v \le 1
$$
is a standard Brownian bridge in terms of which the path of $\br$ on $[u,\hu]$ is represented as
$$
\br(\pu) = x + 
\sqrt{\hu - u }
X^* \left( \frac{\pu - u}{\hu - u} \right)
+
\frac{\pu - u}{\hu - u} ( \hx - x),  \  \pu \in [u, \hu].
$$
On the other hand, with $(t, y) = \mathcal T(u, x)$ and $(\hat t, \hat y) = \mathcal T(\hat u, \hat x)$  the process
$$
Y^*(w):= \frac{\widehat B( t + w ( \hhat - t )) - y - w ( \hy - y )}{\sqrt{\hhat - t }}, ~~~~ 0 \le w \le 1
$$
is another standard Brownian bridge in terms of which the path of $\widehat B$ on $[t,\hhat]$ is represented as
$$
\widehat B(\pt) = y  + \sqrt{\hhat - t } \  Y^* \left( \frac{\pt - t}{\hhat - t} \right)
+
\frac{\pt - t}{\hhat - t} ( \hy - y), ~~~~~~ \pt \in [t, \hhat].
$$
It follows that for an arbitrary choice of $ 0 \le u < \hu < 1$ and real values of $x$ and $\hx$, with $(t,y)$ and $(\hhat, \hy)$ the images of $(u,x)$ and $ (\hu,\hx)$ respectively via the space-time transformation $\mathcal T $ in  \re{spacetime},
\begin{quote}
the image via the space-time transformation of a Brownian bridge from $(u,x)$ to $(\hu,\hx)$ is
a Brownian bridge from $(t,y)$ to $(\hhat,\hy)$. 
\end{quote}
The key observation is that similarly, 
\begin{quote}
the image via the space-time transformation of the
straight line from $(u,x)$ to $(\hu,\hx)$  minus a Brownian excursion of length $\hu - u$ shifted to start at time $u$,
is a straight line from $(t,y)$ to $(\hhat,\hy)$ minus a Brownian excursion of length $\hhat - t$, shifted to start at time $t$.
\end{quote}
Intuitively this is clear from the bridge result, by conditioning each of the bridges to stay above the line joining its endpoints.
This can be made rigorous by a weak convergence argument, conditioning one of the bridges to go no more than a small distance $\epsilon$ above the line, 
passing to the limit as $\epsilon \te 0$, and appealing to the result of \cite{blum_exc83}.  \hfill $\Box$

\medskip

The limiting argument at the end of the previous proof can also be reduced to an invariance of laws of standard Brownian bridges and excursions under a family of space-time transformations indexed by a pair of parameters $0 < u < \hu < 1$.
To see this, observe that the relation between $\br$ and $\widehat B$ implies that the standard bridge $X^*$ derived from $\br$ on $[u,\hu]$
and the standard bridge $Y^*$ derived from $\widehat B$ on $[t,\hhat]$ are related by
\eq
\lb{bridgetrans}
X^*(v) = \frac{ 1 - u - v ( \hu - u ) } {\sqrt{ (1-u) (1 - \hu ) }} \  Y^* \left( \frac{ (1-u) v }{ 1 - u - v (\hu - u) }\right) , ~~~~ 0 < v < 1.
\en
For the sake of completeness, we give the following lemma, which shows that the mapping from one standard bridge or excursion to another standard bridge or excursion which is implicit in Doob's space-time transformation is non-trivial,  but nonetheless easily checked:

\begin{lmm}{lemtrans}
For each fixed choice of $0 \le u < \hu < 1$ if $(Y^*(v), 0 \le v \le 1)$ is a standard Brownian bridge, 
then $(X^*(v), 0 \le v \le 1)$ defined by \re{bridgetrans} has the same distribution as $(Y^*(v), 0 \le v \le 1)$.
Moreover, the same is true with the standard Brownian bridge replaced by a standard Brownian excursion, or by a standard Bessel bridge of any dimension.
\end{lmm}

\proof
For the standard Brownian bridge, the result can be derived as above, or checked more directly by observing that $Y^*$ is evidently a centered
Gaussian process with continuous paths, so it suffices to check that $\IE ( X^*(v) X^*(w) ) = \IE ( Y^*(v) Y^*(w) ) = v (1-w)$ for $0 < v < w < 1$, 
and this is easily done. The result for a Bessel bridge of dimension $d = 1,2, \ldots$ follows easily from the representation of the square of this
process as the sum of squares of $d$ independent standard Brownian bridges. For $d = 3$ this gives the result for standard excursion, by the well known
identification of the standard excursion with a three-dimensional Bessel bridge due to David Williams.
\endpf



\section{The distribution of $M$}
\label{distM}

\paragraph{Proof of Theorem \ref{Rep1}.} 
Proposition \ref{indep} implies that the distribution of the maximal difference $M$ is the same for 
a standard Brownian motion on $[0,1]$ as for a Brownian bridge from $(0,0)$ to $(1, b)$ for an arbitrary real number $b$. 
So it is enough to establish the characterization of $M$ provided in Theorem \ref{Rep1} with a standard Brownian motion $B$ on $[0,1]$ replaced by a standard Brownian bridge $\br$.

Let $(V_i, i \in \mathbb Z)$ denote the sequence of times of vertices of $\ovC_{[0,1]} \br$ indexed in the same way as in Proposition \ref{indep} and Theorem \ref{conc}, and let $(T_i, i \in \mathbb Z)$ be the corresponding lengths of the segments of $\ovC_{[0,1]} \br$; i.e, $T_i = V_i - V_{i-1}, i \in \mathbb Z$ so the sequence $(T_i, i \in \mathbb Z )$ is independent of the sequence $(M_{3,i}, i \in \mathbb Z )$. 

From the definition of $M$ and Theorem \ref{conc}, we have readily that 
\begin{eqnarray*}
M = \max_{i} \sqrt{T_i} M_{3,i}
\end{eqnarray*}
where the $M_{3,i}$'s are independent random variables identically distributed as $M_3$, the maximum of a standard Brownian excursion of length 1. 

Now let $T^{(1)} \ge T^{(2)} \ge \cdots$ denote the values of $T_i$ put in decreasing order. 
According to the result of Suidan \cite{suidan01}, the distribution of the sequence ($T^{(j)}, j = 1,2 \ldots)$ 
is the limiting distribution of ranked lengths of cycles of a uniform random permutation of $n$ elements, with cycle lengths normalized by $n$,
commonly known as the Poisson-Dirichlet distribution with parameters 0 and $1$. See also \cite{MR994088}, \cite{MR1910531}, \cite{py95pd2}. Suidan gives this result for the concave majorant of Brownian motion on $[0,1]$, but it applies just as well to the Brownian bridge,
by Proposition 4. It is also known \cite{vs77} that 
this asymptotic distribution of cycle lengths is obtained by ranking the terms of a uniform stick-breaking process in decreasing order.
Theorem \ref{Rep1} and Corollary \ref{uniform} follow immediately.
\endpf

\begin{figure}[!h]
\centering
\centerline{\epsfig{file = 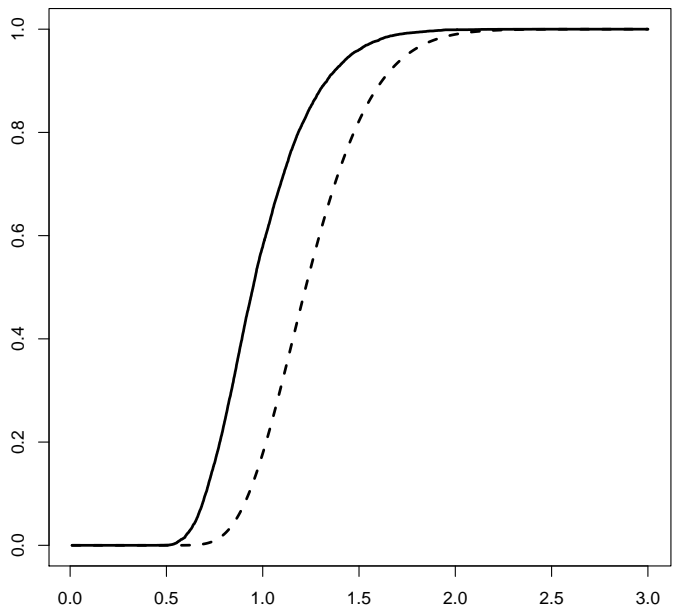, width = 10cm}}
\vspace*{-0.5cm}
\caption{In solid line is the plot of a Monte Carlo approximation of $F_M$ based on a sample of size 5,000. In dashed line is the plot of the distribution function $F_3$. As the figure suggests, $F_3(u) \le F_M(u)$  for all $u$. This is an immediate consequence of (\ref{m3diff}) and the definition of $M$.}
\label{FM}
\end{figure}

\begin{figure}[!h]
\centering
\centerline{\epsfig{file = 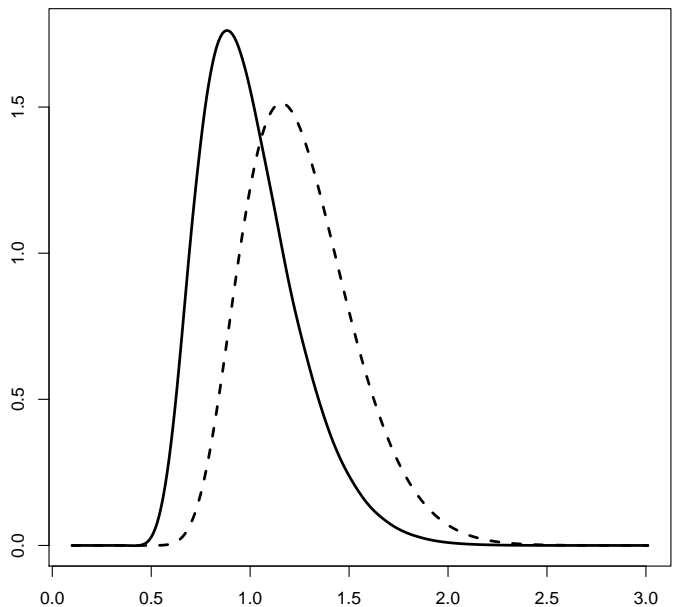, width = 10cm}}
\vspace*{-0.5cm}
\caption{In solid line is the plot of a Monte Carlo approximation of $f_M$ based on a sample of size 5,000. In dashed line is the plot of the density function $f_3$.}
\label{fM}
\end{figure}

Theorem \ref{Rep1} offers an easy way of simulating values of $M$. For the Monte Carlo implementation, we can use
the representation \re{m3diff} of $M_3$
and the Donsker approximation of $\br$.  It follows also from Theorem \ref{Rep1} that the distribution and density functions of $M$, denoted hereafter by $F_M$ and $f_M$ respectively, are given by
\begin{eqnarray}\label{ExpFM}
F_M(x) = E\left[\prod_i F_3\left( \frac{x}{\sqrt{L_i}} \right)  \right],  ~~~~~  x > 0
\end{eqnarray}
and
\begin{eqnarray}\label{ExpfM}
f_M(x) = \sum_{i} E\left[\frac{1}{\sqrt{L_i}} f_3\left( \frac{x}{\sqrt{L_i}} \right) \prod_{j \ne i}  F_3\left( \frac{x}{\sqrt{L_j}} \right) \right],  ~~~~~  x > 0
\end{eqnarray}
where $F_3$ and $f_3$ are the distribution and density functions of $M_3$, known to be given by
\begin{eqnarray}\label{DisM3}
F_3(y)  = 1 - 2 \sum_{n =1}^\infty ( 4 n^2 y^2 - 1) e^{- 2 n^2 y^2} ,  ~~~~~  y > 0
\end{eqnarray}
and
\begin{eqnarray*}\label{DensM3}
f_3(y)  = 8\sum_{n = 1}^\infty n^2 y (4n^2 y^2 - 3)  e^{- 2 n^2 y^2} ,  ~~~~~  y > 0,
\end{eqnarray*}
see e.g. \cite{MR0402955}, \cite{bpy99z} for further information about this distribution which has been extensively studied.  
Monte Carlo approximations of the distribution function $F_M$ and its density $f_M$ based on 5,000 independent copies of the uniform stick breaking process are shown in Figure \ref{FM} and Figure \ref{fM} along with the distribution and density functions of $M_3$. 

\bigskip

Another technique for getting at features of the distribution of $M$ is to introduce  a standard exponential variable $\gamma_1$ independent of $M$, and observe that 
\eq
\lb{maxdiff1x}
\sqrt{\gamma_1} M \ed \max_i \sqrt{\gamma_1 L_i} M_{3,i}
\en
where the $\gamma_1 L_i$ are the points of a Poisson point process with intensity $x^{-1} e^{-x} dx$ for $x >0$,
independent of the $M_{3,i}$'s. These are the jumps of a gamma process $(\gamma_s, 0 \le s \le 1)$ from which the $L_i$ can be recovered
by first normalizing the jumps by $\gamma_1$ and then putting the jumps in size-biased random order.
This leads to the following proposition.

\begin{prp}{gamma1}
For all $x > 0$
\begin{eqnarray*}
\lb{maxdiff1}
\IP ( \sqrt{\gamma_1} M  \le x) = \exp( - \nu(x,\infty) )
\end{eqnarray*}
where
\eq
\lb{maxdiff12}
\nu(x,\infty) = \int_0 ^\infty y^{-1} e^{-y} \ \left[1 - F_3\left( \frac{x}{\sqrt{y}} \right)\right] dy.
\en
\end{prp}

\medskip

\par \noindent \textbf{Proof.}  \ Note that the event $\{ \sqrt{\gamma_1} M  \le x \}$ occurs if and only if there is no point $\sqrt{\gamma_1 L_i} M_{3,i}$ in the interval $(x,\infty)$. Since these $\sqrt{\gamma_1 L_i} M_{3,i}$ are the points of a Poisson point process with intensity measure $\nu$; the image via the map $(x,m) \mapsto \sqrt{x} m $ of the measure 
$$x^{-1} e^{-x} \IP (M_3 \in dm),$$  
it follows that
\begin{eqnarray*}
\IP ( \sqrt{\gamma_1} M  \le x) = \exp( - \nu(x,\infty) )
\end{eqnarray*}
with $\nu(x, \infty)$ given in \re{maxdiff12}. \hfill $\Box$ 

\medskip

See also \cite{ap02d,ap01d,py97rh} for similar calculations.

\bigskip

The following theorem gives another representation of $F_M$ and $f_M$. Let $K_m$ be the modified Bessel function of the second kind and order $m \in [0, \infty)$, and consider the function 
\begin{eqnarray}\label{Bess}
G(t) = \prod_{n=1}^\infty \exp\left\{-4 \left[ 2 \sqrt{2} t \ n K_1( 2 \sqrt{2} t \ n) - K_0( 2 \sqrt{2} t \ n) \right] \right\}, ~~~t > 0.
\end{eqnarray}
Also, let $\mathcal L^{-1}$ be the operator of inverse Laplace transform.

\begin{thm}{LapInv} \
\begin{description}
\item (i) For any real $r > 0$, we have 
\begin{eqnarray}\label{IdMoments}
E(M^r) = \frac{2}{\Gamma(r/2)} \int_0^\infty t^{r-1} (1 - G(t)) dt.
\end{eqnarray}

\item (ii) For all $x > 0$ the distribution and density functions $F_M$ and $f_M$ are given  respectively by
\begin{eqnarray}\label{lapinvDis}
F_M(x) = \left[\mathcal{L}^{-1}(G(\sqrt{t})/t)\right]\left(\frac{1}{x^2}\right)
\end{eqnarray}
and
\begin{eqnarray}\label{lapinvDens}
f_M(x) = \frac{2}{x^3}\left[\mathcal{L}^{-1}(1-G(\sqrt{t}))\right]\left(\frac{1}{x^2}\right).
\end{eqnarray}
\end{description}
\end{thm}

\medskip

\par \noindent \textbf{Proof.} \ It follows from \re{DisM3} and \re{maxdiff12} that 
\begin{eqnarray*}
\nu(x, + \infty) & = &  2 \sum_{n=1}^\infty \int_0^\infty y^{-1} e^{-y} (4n^2 x^2 y^{-1} - 1) e^{-2 n^2 x^2 y^{-1}} dy   \\
                 & = &  2 \sum_{n=1}^\infty \int_0^\infty t^{-1} e^{-2 n^2 x^2 t} (2 t^{-1} - 1) e^{- t^{-1}} dt \\
                 & = &  2 \sum_{n=1}^\infty (2 A_n(x) - B_n(x))   \\
\end{eqnarray*}
where
\begin{eqnarray}\label{An}
A_n(x) = \left[\mathcal{L}(e^{-t^{-1}}/t^2)\right](2 n^2 x^2) = 2\sqrt{2} n x K_1(2 \sqrt{2} n x) 
\end{eqnarray}
and 
\begin{eqnarray}\label{Bn}
B_n(x) = \left[\mathcal{L}(e^{-t^{-1}}/t)\right](2 n^2 x^2) = 2 K_0(2\sqrt{2} n x),
\end{eqnarray}
so that
$$
\exp(- \nu(x, \infty)) = \prod_{n=1}^\infty \exp\left\{-4 \left[ 2 \sqrt{2} x \ n K_1( 2 \sqrt{2} x \ n) - K_0( 2 \sqrt{2} x \ n) \right] \right\} = G(x).
$$
The derivation of the expressions in (\ref{An}) and (\ref{Bn}) is deferred to the appendix.  Now for $r> 0$
\begin{eqnarray*}
E(\gamma^{r/2}_1 M^r)& =&  \int_0^\infty \left(1 - F_{\gamma^{r/2}_1 M^r}(x)\right) dx \\
                     & = & \int_0^\infty \left(1 - \exp(\nu(x^{1/r}, \infty)) \right) dx \\
                     & =  & \int_0^\infty \left(1 - G(x^{1/r}) \right) dx = \int_0^\infty r x^{r-1}(1 - G(x)) dx.  
\end{eqnarray*}
The claim in (i) follows now from independence of the random variables $\gamma_1$ and $M$.

To show (ii), we use again the independence of $\gamma_1$ and $M$. We can write 
\begin{eqnarray*}
\IP( \sqrt{\gamma_1} M  \le x) = \int_0^\infty F_M\left( \frac{x}{\sqrt{s}} \right) e^{-s} ds.
\end{eqnarray*}
Using the change of variable $t = s/\sqrt{x}$, we get
\begin{eqnarray*}
\int_0^\infty F_M\left( \frac{1}{\sqrt{t}} \right) e^{-t x} dt =  \exp(-\nu(\sqrt{x}, \infty))/x.
\end{eqnarray*}
Thus, for all $t > 0$
$$
F_M\left( \frac{1}{\sqrt{t}} \right) = \left[\mathcal{L}^{-1}\left( \exp(-\nu(\sqrt{x}, \infty))/x \right)\right](t)
$$
and the expression of $F_M$ follows. The expression $f_M$ can be obtained immediately by using known properties of the operator $\mathcal{L}$. \hfill $\Box$

\begin{table}[!h]
\caption{Approximation of the moments of $M$ of order $r=1, ..., 8$.  }
\begin{center}
\begin{tabular}{ccc}
\hline
$r$  & $E(M^r)$  &  $\overline{X}^r$\\
\hline
1 & 0.999399  &  0.997366 \\
2 & 1.060258 &  1.056803\\
3 & 1.195155 &  1.190869\\
4 & 1.431334 & 1.427101\\
5 & 1.819154 & 1.816777 \\
6 & 2.448679 & 2.452149\\
7 & 3.481508  & 3.499897\\
8  &5.212503  &  5.266828\\
\hline
\end{tabular}
\end{center}
\label{Moments}
\end{table}

\medskip
\medskip

The identity in (\ref{IdMoments}) gives a way of calculating the $r$th moment of $M$ via numerical integration. We used this approach to compute $E(M^r)$ for $r=1, .., 8$, and the values are reported in the first column of Table \ref{Moments}. In the second column are the corresponding empirical estimators for the same moments using 20,000 simulated independent copies of $M$. On the other hand, the inverse Laplace transforms in the expressions of $F_M$ and $f_M$ given in (\ref{lapinvDis}) and (\ref{lapinvDens}) can be approximated with very high precision using for example the Gaver-Stehfest algorithm. We refer to \cite{MR2274708} for a detailed description of this method as well as a nice summary of other classical algorithms used to approximate inverse of Laplace transforms. However, implementation of the Gaver-Stehfest algorithm  requires using a multiple precision software or some adequate arbitrary precision library. This will be pursued elsewhere.

\section{Complements}
\label{complements}
Groeneboom's description of $\ovC_{[0, \infty)} B$ implies that for each fixed $t > 0$ the joint density 
of $V^-_t$, the time of the last vertex of $\ovC_{[0, \infty)} B$ before time $t$, and
$V^+_t$, the time of first vertex of $ \ovC_{[0, \infty)} B $ after time $t$ is given by the formula
\begin{eqnarray*}
f_{V^-_t, V^+_t}(v_1, v_2) = \frac{2}{(v_2 - v_1)^{3/2}} E\left[Z_+ \left(\frac{X}{\sqrt{v_1}} - \frac{Z}{\sqrt{v_2 - v_1}}  \right)_+\right] 1_{\{0 < v_1 < t < v_2\}}
\end{eqnarray*}
where $y_+ = y 1_{y \ge 0}$ and $X$ and $Z$ are independent standard normal variables. 
We show in an appendix (Proposition \ref{Explicit}) that this joint density can be presented more explicitly as
\begin{eqnarray}\label{straddt}
f_{V^-_t, V^+_t}(v_1, v_2) = \frac{1}{\pi (v_2 - v_1)^{2} } \left( \sqrt{\frac{v_2 - v_1}{v_1}} - \arctan\left(  \sqrt{\frac{v_2 - v_1}{v_1}}\right) \right) 1_{\{0 < v_1 < t < v_2\}}.
\end{eqnarray}

Using the time transformation $t \mapsto t/(t+1)$, the joint density of $X_u$ and $Y_u$, the last and first times vertices of the concave majorant of $\br$ occurring before and after $u \in (0,1)$, is given by
\begin{eqnarray}\label{straddu}
f_{X_u, Y_u}(x, y) &= & \frac{1}{\pi (y-x)^2} \left(\sqrt{\frac{y-x}{x(1-y)}} - \arctan\left(\sqrt{\frac{y-x}{x(1-y)}}\right)\right) \nonumber \\ && \hspace{0.2cm} \times \ 1_{\{0 < x < u < y < 1\}}.
\end{eqnarray}

Now let $L_1$ be the length of the segment of the \cmaj \ of a standard Brownian bridge $\br$ on $[0,1]$ covering a uniform random number $U_1$ independent of $\br$, and let us verify using this formula that the distribution of $L_1$
is also uniform on $[0,1]$, as shown already by Corollary \ref{uniform}.
Using the expression in (\ref{straddu}), we find after some algebra that $L_1$ has density
\begin{eqnarray*}
f_{L_1}(l) &= & \int_0^{1-l} \frac{1} {\pi l} \left(\sqrt {\frac{l}{x (1- l - x)}}   - \arctan\left( \sqrt {\frac{l}{x (1- l - x)}} \right)  \right) dx \\
& = & \frac{1-l} {\pi l}  \int_0^{1} \left(k \frac{1}{\sqrt {l (1- l)}}   - \arctan\left( \frac{k}{\sqrt{u (1-u)}} \right)  \right) du,
\end{eqnarray*}
using the change of variable \ $u = x/(1-l)$ and putting $k = \sqrt{l}/(1-l)$.

Now,
\begin{eqnarray*}
\int_0^{1} \frac{du}{\sqrt {u (1- u)}} = \int_0^{\pi/2} 2 d\theta = \pi, \ \textrm{ by the change of variable \ $u = \sin^2(\theta) $},
\end{eqnarray*}
Using the same change of variable we can write
\begin{eqnarray*}
\int_0^{1}\arctan\left( \frac{k}{\sqrt{u (1-u)}} \right)  du &= & \int_0^{\pi/2}\arctan\left( \frac{2k}{\sin(2\theta)}\right) \sin(2\theta) d\theta\\
& = & \frac{1}{2} \int_0^{\pi}\arctan\left( \frac{2k}{\sin(t)}\right) \sin(t) dt\\
& = & \frac{1}{2} \left[ - \cos(t) \arctan\left(\frac{b}{\sin(t)}\right) \right]_0^\pi \\
&&  - \frac{1}{2} \  b \int_0^{\pi} \frac{\cos^2(t)}{\sin^2(t) + b^2 } dt, \  \textrm{with \ $b = 2k$} \\
&= & \frac{\pi}{2} -  b \int_0^{\pi/2}  \frac{\cos^2(t)}{\sin^2(t) + b^2 } dt \\
& =  & \frac{\pi}{2} -  \frac{b}{b^2 + 1}  \int_0^{\pi/2}  \frac{dt}{\tan^2(t) + b^2/(b^2 +1) }
\end{eqnarray*}
where the last equality follows from the identity $1/\cos^2(t) = 1 + \tan^2(t)$. Put $c =b^2/(b^2 +1) $. Using the change of variable $\tan(t) = x$, we get
\begin{eqnarray*}
\int_0^{\pi/2}  \frac{dt}{\tan^2(t) + c} &  = & \int_0^\infty \frac{dx}{(1+x^2)(x^2 + c)} =  \frac{1}{1 - c} \int_0^\infty \left(\frac{1}{x^2+ c} - \frac{1}{1+x^2} \right) dx \\
 & = & \frac{1}{1 - c} \frac{\pi}{2} \left( \frac{1}{\sqrt{c}} - 1 \right).  
\end{eqnarray*}
Now,
\begin{eqnarray*}
c = \frac{4l}{2l + 1 + l^2} = \frac{4l}{(l+1)^2} 
\end{eqnarray*}
and hence
\begin{eqnarray*}
\int_0^{1}\arctan\left( \frac{k}{\sqrt{u (1-u)}} \right)  du & = &\frac{\pi}{2} \left[1 - \frac{2 \sqrt{l}}{1-l} \left( \frac{l+1}{2\sqrt{l}} - 1 \right) \right]
 =  \frac{\pi \sqrt{l}}{1 + \sqrt{l}}
\end{eqnarray*}
which implies that $f_{L_1}(l) = 1$ for $0 < l < 1$. Thus we have an independent verification  of the uniform
distribution of $L_1$ asserted by Corollary \ref{uniform}. But it is far from clear how to derive the further
results of Corollary \ref{uniform} by this approach.

\paragraph{Acknowledgments.} The first author would like to thank H. Doss for interesting discussions.

\section{Appendix}


We collect in this appendix some further identities and computations of an analytic kind which arise in connection
with the material of the paper.

\begin{prp}{Explicit2} 
Let $X$ and $Y$ have the standard bivariate normal distribution with correlation $E(XY) = \rho$. Then
\eq
\label{idrosen}
E( X_+ Y_+ ) = \frac{1}{2\pi} \sqrt{1 - \rho^2} + \rho P(X > 0 , Y > 0 ) 
\en
where
\begin{eqnarray}
\label{id2}
P(X > 0 , Y > 0 )  =  \frac{1}{2\pi} \left( \frac{\pi}{2} + \arctan (\rho/\sqrt{1 - \rho^2}) \right) 
\end{eqnarray}
which for $\rho < 0$ can be rewritten as
\eq
\label{id2neg}
P(X > 0 , Y > 0 )  =  -\frac{1}{2\pi} \arctan (\sqrt{1 - \rho^2}/\rho).
\en
\end{prp}
\proof
Rosenbaum \cite[formula (5)]{rosenbaum61} gives an expression for
$E ( X Y 1_{\{X >h, Y > k \}})$ in terms of the probability $P (X >h, Y > k)$ and the standard normal  density and distribution
functions. In the present case $h = k = 0$ and Rosenbaum's formula simplifies  to \re{idrosen}.
Formula \re{id2} is well known, see e.g. \cite{stuartord}.
If $\rho <0$ the passage from \re{id2} to \re{id2neg} is made via the 
trigonometric identity $\arctan(x) + \arctan(1/x) = - \pi/2, \ \ \forall \ x < 0$. \hfill $\Box$

\medskip

\begin{prp}{Explicit} Let $a, b > 0$, and $Z$ and $W $ independent standard normal variables. Then
\begin{eqnarray*}
E\left[Z_+ \left(\frac{W}{a} - \frac{Z}{b}\right)_+    \right] = 
\frac{1}{2 b \pi} \left(\frac{b}{a} - \arctan\left(  \frac{b}{a}\right)\right)
\end{eqnarray*}
where $y_+ = y 1_{y \ge 0}$.
\end{prp}
 
\medskip

\par \noindent \textbf{Proof.} \ Let $t:= b/a$, and $X: = Z$ and $Y = (t W  - Z)/\sqrt{1+ t^2}$.  Since  $(X, Y)$ is a  standard bivariate normal with correlation $ \rho = -1/\sqrt{ t^2 +1} < 0 $ where $t =  - \sqrt{ 1 - \rho ^2} / \rho > 0$  then it follows from Proposition \ref{Explicit2} that 
\begin{eqnarray*}
E\left[Z_+ \left(t W  - Z\right)_+    \right] & = &  -\frac{1}{2\pi \rho} \left( \sqrt{1 - \rho^2} - \rho \arctan(\sqrt{1 - \rho^2}/\rho) \right)
  \\
& = & \frac{1}{2 \pi} \left( t - \arctan (t) \right)  
\end{eqnarray*}
and the conclusion follows after replacing $t$ by $b/a$.   \hfill $\Box$

\medskip

\paragraph{Proof of the identities (\ref{An}) and (\ref{Bn}).} \ We start by showing (\ref{Bn}). Letting $ z= 2 \sqrt{2} n x$, we need to show that 
\begin{eqnarray*}
K_0(z) = \frac{1}{2} \int_0^\infty \frac{e^{-1/t}}{t} \ e^{-z^2 t/4 } dt.
\end{eqnarray*} 
We have that 
\begin{eqnarray*}
\int_0^\infty \frac{e^{-1/t}}{t} \ e^{-z^2 t/4 } dt = e^{-z} \int_0^\infty \frac{1}{t} \ e^{-\frac{(zt-2)^2}{4t}} dt  
\end{eqnarray*} 
and putting  $u = (zt - 2)/(2 \sqrt{t})$, it follows that $t = (u + \sqrt{u^2 + 2z})^2/z^2$ and $dt = 2t du/\sqrt{u^2 + 2z}$. Hence
\begin{eqnarray*}
\int_0^\infty \frac{e^{-1/t}}{t} \ e^{-z^2 t/4 } dt = 2 e^{-z} \int_{-\infty}^\infty \frac{e^{-u^2}}{\sqrt{u^2 +  2z}} du.  
\end{eqnarray*} 
It follows that 
\begin{eqnarray*}
\frac{1}{2} \int_0^\infty \frac{e^{-1/t}}{t} \ e^{-z^2 t/4 } dt = e^{-z} \int_{-\infty}^\infty \frac{e^{-u^2}}{\sqrt{u^2 +  2z}} du = K_0(z)
\end{eqnarray*} 
see e.g. Hunter \cite{MR0158104}.  To show (\ref{An}), it is enough to show that
\begin{eqnarray*}
\int_0^\infty \frac{e^{-1/t}}{t^2} \ e^{-z^2 t/4 } dt = z K_1(z).
\end{eqnarray*}
Using the fact $\lim_{z \to 0} z K_1(z) = 1$ (see e.g. Mechel \cite{MR0202282}), and noting that the integral on the left hand side is equal to 1 for $z=1$ it is enough to show that
\begin{eqnarray*}
- \frac{z}{2} \int_0^\infty \frac{e^{-1/t}}{t} \ e^{-z^2 t/4 } dt = (z K_1(z))'.
\end{eqnarray*}
The calculations above imply that 
\begin{eqnarray*}
- \frac{z}{2} \int_0^\infty \frac{e^{-1/t}}{t} \ e^{-z^2 t/4 } dt = - z K_0(z).  
\end{eqnarray*}
We conclude by using the well-known identity $(zK_n(z)' =  -z K_{n-1}(z)$ where $K_m$ is the modified Bessel function of the second kind of order $m$.  \hfill $\Box$

\medskip

\end{document}